\documentclass{article}
\usepackage{amssymb,amsmath,amsthm}
\newtheorem{theorem}{Theorem}
\newtheorem{definition}{Definition}
\newtheorem{lemma}{Lemma}
\newtheorem{proposition}{Proposition}
\newtheorem{corollary}{Corollary}
\newtheorem{remark}{Remark}
\date{}
\numberwithin{equation}{section} \numberwithin{theorem}{section}
\numberwithin{lemma}{section} \numberwithin{corollary}{section}
\numberwithin{remark}{section} \numberwithin{proposition}{section}
\numberwithin{definition}{section}
\begin{document}
\newcommand{\n}{\noindent}
\newcommand{\vs}{\vskip}
\title{Porosity of the Free Boundary in a Class of Higher-Dimensional Elliptic Problems}

\vs 0.5cm
\author{Abdeslem Lyaghfouri \\
United Arab Emirates University\\
Department of Mathematical Sciences\\
Al Ain, Abu Dhabi, UAE} \maketitle

\vskip 0.5cm
\begin{abstract}
We investigate a class of $n$-dimensional ($n\geq 2$) free boundary elliptic 
problems which includes the dam problem, the aluminum problem, and the lubrication problem.
We establish that the free boundary in this class is a porous set, which implies
its Hausdorff dimension being less than $n$, which in turn leads to
its Lebesgue measure being zero. Our proof relies on the comparison 
of the solution with an appropriately constructed barrier function.
\end{abstract}
\vs 0.5cm \noindent AMS Mathematics Subject Classification:
35J15, 35R35

\vs 0.5cm
\section{Introduction}\label{1}

\vskip 0.5cm

Let $\Omega$ be the open set of $\mathbb{R}^n$ defined by
$$\Omega=\{(x',x_n)\in \mathbb{R}^n\,\,/\,\,x'=(x_1,...,x_{n-1})\in B'_\rho,\,\, \gamma(x')<x_n<\gamma(x')+l\}$$

\n and let $\Gamma=\{(x',\gamma(x')+l)\,\,/\,\,x'\in B'_\rho\}$, where $\rho$, $l$ are two positive numbers, 
$B_\rho'=\{x'\in \mathbb{R}^{n-1}\,\,/\,\,|x'|<\rho\}$, and $\gamma:~B_\rho'~\rightarrow~\mathbb{R}$ is a $C^{1,1}$ function. 

\vs0.2cm\n Let $a(x)= (a_{ij}(x))$ be an $n\times n$ matrix and $h(x)$ denotes a scalar function such that 
the following conditions hold for positive constants $\lambda$, $\Lambda$, $\bar h$, $\underline{h}$ and $p>n$:
\begin{eqnarray}\label{e1.1-1.5}
 & \sum_{i,j}|a_{ij}(x)|\leq \Lambda,\quad\text{  for a.e. }x\in \Omega,\\
&{ a}(x)\xi.\xi \geq \lambda\vert \xi\vert^2,\quad\text{ for all }\xi\in
\mathbb{R}^{n}\quad\text{and for a.e. }x\in \Omega,\\
& \underline{h} \leq h(x)\leq \bar h \quad\text{  for a.e. }x\in \Omega  \\
& h_{x_n}\in L_{loc}^p(\Omega\cup\Gamma).\\
& h_{x_n}(x)\geq 0\quad\text{  for a.e. }x\in \Omega  .
\end{eqnarray}

\n Let $e_n=(0,...,0,1)$. We are interested in investigating the following problem:
\begin{equation*}(P)
\begin{cases}
& \text{Find } (u, \chi) \in  H^{1}(\Omega)\times L^\infty (\Omega)
\text{ such that } :\\
& (i)\quad u=0\quad\text{on}\quad \Gamma\\
& (ii)\quad u\geq 0, \quad 0\leq \chi\leq 1 ,
\quad u(1-\chi) = 0 \,\,\text{ a.e.  in } \Omega\\
 & (iii)\quad \displaystyle{\int_\Omega }\big( a(x) \nabla
u + \chi h(x)e_n\big) .\nabla\xi dx \,\leq\, 0\\
&\hskip 1.7cm \forall \xi \in H^{1}(\Omega),\quad\xi = 0 \text{ on }
\partial\Omega\setminus\Gamma, \quad\xi\geq 0 \text{ on } \Gamma.
\end{cases}
\end{equation*}

\n This formulation encompasses various free boundary problems, notably the dam problem 
(references \cite{[A1]}, \cite{[AG]}, \cite{[CC]}, \cite{[CL2]}, \cite{[L1]}, \cite{[SV]}), where 
$\Omega$ denotes a porous medium, $a(x)$ represents its permeability matrix, and $h(x)=a_{nn}(x)$.
Additionally, the weak formulation extends to include the lubrication problem (reference \cite{[AC]}), 
with $a(x)=h^3(x)I_2$, where $I_2$ denotes the $2\times2$ identity matrix, and $h(x)$ is a scalar 
function associated with the Reynolds equation. A third model addressed herein is the aluminum 
electrolysis problem (references \cite{[BMQ]}, \cite{[CM]}), in which $a(x)$ is a variable
diagonal matrix.

\vs 0.2cm\n Our primary objective in this study is to investigate the free boundary $(\partial\{u>0\})\cap\Omega$, 
denoted in this paper by $FB$, which delineates the interface between the fluid-containing 
region and the fluid-free region in scenarios such as dam and lubrication problems, 
as well as in the context of aluminum electrolysis, where it separates 
regions with liquid and solid aluminum. The regularity of the set $FB$ has been explored by numerous 
researchers in various contexts. 
For instance, in \cite{[A2]}, H.W. Alt established that it is a Lipschitz continuous curve 
$x_n=\Phi(x')$ under the condition that $a(x)$ is the identity matrix.
In the two dimensional case, A. Lyaghfouri \cite{[L1]} demonstrated that $FB$ is the graph of a continuous 
curve $x_2=\Phi(x_1)$, provided that $a_{12}(x)=0$ and $h(x)=a_{22}(x)$ is nondecreasing with respect to $x_2$.

\vs 0.2cm\n When $a$ and $h$ are both non-constant and unrelated to each other, proving the continuity
of the free boundary becomes more challenging. This complexity was initially addressed in the two-dimensional 
context by M. Chipot \cite{[C]}, and later explored by S. Challal and the author \cite{[CL1]} under varied assumptions. 
Moreover, under general assumptions, the authors demonstrated in \cite{[CL3]}, that $FB$ is represented locally by a family of 
continuous functions. However, all these proofs primarily focused on two dimensions and were not 
easily extendable to higher dimensions. Given that under assumption (1.5), $FB$ is the graph of a lower semi-continuous
function $\phi(x')$, a pivotal aspect of demonstrating the continuity, was to establish the upper semi-continuity 
of $\phi$, which in dimension two rested on the ability to identify two points in arbitrarily close proximity to a free boundary 
point, such that $u$ vanishes above these points. This facilitated the construction of a barrier function for 
comparison with $u$ on a rectangle positioned above the free boundary point. However, this approach becomes 
ineffectual in dimensions $n\geq 3$.

\vs0.2cm\n In this paper, our aim is to investigate the regularity of the set $FB$ within the broader framework 
of $n$-dimensional case ($n \geq 2$). Our main contribution lies in establishing the porosity 
of $FB$, a result that notably entails its Hausdorff dimension being less than $n$, consequently resulting 
in its Lebesgue measure being zero.

\vs 0.5cm\n We begin with two remarks.

\begin{remark}\label{r1.1}
Consider the change of variables $y=\Phi(x)=(\Phi_i(x))$,
$x=\Psi(y)=(\Psi_i(y))$, where

\[\Phi_i(x)=\left\{
\begin{array}{ll}
    x_i, \hbox{ if }~~ i=1,...,n-1&\\
    x_n-\gamma(x_1,...,x_{n-1}), \hbox{ if }~~ i=n,&  \\
\end{array}%
\right.\] \n and
\[\Psi_i(y)=\left\{
\begin{array}{ll}
    y_i, \hbox{ if }~~ i=1,...,n-1&\\
    y_n+\gamma(y_1,...,y_{n-1}),\hbox{ if }~~ i=n. &  \\
\end{array}%
\right.\]

\n Then $\Psi:\,U=B_\rho'\times(0,l)\rightarrow \Omega\, $ and
$\Phi:\,\Omega\rightarrow U\,$ are $C^{1,1}$ diffeomorphisms with $\Psi=\Phi^{-1}$ and $J_\Phi(x)=J_\Psi(y)=1$ for all
$x\in\Omega$ and all $y\in U$. Moreover, setting $b(y)=D\Phi(\Psi(y))a(\Psi(y))^tD\Phi(\Psi(y))$, $v(y)=uo\Psi(y)$,
$\theta(y)=\chi o\Psi(y)$, and $k(y)=ho\Psi(y)$, it is not difficult to verify that $(v,\theta)$ is
a solution of the following problem:

\begin{equation*}(Q)
\begin{cases}
& \text{Find } (v, \theta) \in  H^{1}(U)\times L^\infty (U)
\text{ such that}:\\
& (i)\quad v=0\quad\text{on}\quad S=B_\rho'\times\{l\}\\
& (ii)\quad  v\geq 0, \quad 0\leq \theta\leq 1 ,
\quad v(1-\theta) = 0 \,\,\text{ a.e.  in } U\\
 & (iii)\quad \displaystyle{\int_U }\big( b(y) \nabla
v + \theta k(y)e_n\big) .\nabla\zeta dy \,\leq\, 0\\
&\hskip 1.7cm \forall \zeta \in H^{1}(U),\quad\zeta = 0 \text{ on }
\partial U\setminus S, \quad\zeta\geq 0 \text{ on } S.
\end{cases}
\end{equation*}

\n Therefore, without sacrificing generality, we may consider $\Omega$ to be 
a cylinder $B_\rho'\times(0,l)$, and $\Gamma$ the flat surface $B_\rho' \times {l}$. 
This assumption will persist throughout the paper. 
\end{remark}

\begin{remark}\label{r1.2} $i)$ Henceforth, we denote by $(u,\chi)$ a solution to problem $(P)$

\vs 0.2cm \n $ii)$ From $(P)iii)$, we derive the equation
$div\big(a(x) \nabla u\big)=-(h\chi)_{x_n}$ satisfied in distributional sense.
Since $\chi \in L^\infty(\Omega)$ and due
to (1.1)-(1.3), it follows (see \cite{[GT]} Theorem 8.24, p. 202)
that $u$ is locally H\"{o}der continuous in $\Omega$.

\vs 0.2cm \n $iii)$ Since $u=0$ on $\Gamma$, we also have (see \cite{[GT]} Theorem 8.29, p. 205)
that $u$ is locally H\"{o}der continuous in $\Omega\cup\Gamma$.
\end{remark}

\vs 0,3cm\n Prior to presenting our principal contribution, it is essential to recall a definition (see \cite{[Z]} for example ).

\vs0.3cm
\begin{definition}\label{d1.1}
We define a set $E\subset \mathbb{R}^n$ as porous with a porosity parameter
$\delta>0$ if there exists a positive number $r_0$ such that for every point
$x\in E$ and every radius $r\in (0,r_0)$, there exists a point $y\in\mathbb{R}^n$
such that the open ball $B_{\delta r}(y)$ is entirely contained within the complement of
$E$ in the ball $B_{ r}(x)$, i.e., $B_{\delta r}(y)\subset
B_{ r}(x)\setminus E$.
\end{definition}

\vs0.2cm

\begin{remark}\label{r1.3} A porous set with porosity $\delta$ possesses a Hausdorff dimension that does not surpass 
$n-c\delta^n$ (see \cite{[MV]}, \cite{[T]}), where $c=c(n)>$ represents a constant dependent 
solely on $n$. Consequently, a porous set has Lebesgue measure zero.
\end{remark}

\vs 0,3cm\n We shall additionally presume that

\begin{equation}\label{e1.6}
a\in C^{0,\alpha}_{loc}(\Omega\cup\Gamma)\quad\text{for some }\alpha\in(0,1).
\end{equation}

\vs 0,3cm\n Under the aforementioned assumptions, our main result is presented below.

\begin{theorem}\label{t1.1} The free boundary $(FB)$ is locally a porous set.
\end{theorem}
 
\begin{corollary}\label{c1.1}
\[\chi=\chi_{\{u>0\}}\quad\text{a.e. in }  \Omega\]
\end{corollary}

\vs 0,3cm \n We start by recalling a crucial regularity result: the Lipschitz continuity of
$u$. Under assumptions (1.1)-(1.4) and (1.6), we have
(references \cite{[CL1]}, \cite{[L2]} and \cite{[L3]}):

\begin{proposition}\label{p1.1}
\[u\in C^{0,1}_{loc}(\Omega\cup\Gamma)\]
\end{proposition}

\vs 0.3cm

\begin{remark}\label{r1.4}
The Lipschitz continuity of $u$ stands as a pivotal element in establishing our main result. It's noteworthy that this level of regularity is optimal, as it is directly tied to the discontinuity exhibited by $\nabla u$ across the set $FB$. Additionally, it's vital to highlight that the Lipschitz continuity of $u$ has played a key role in proving the regularity results concerning the free boundary in both \cite{[CL1]} and \cite{[CL3]}.
\end{remark}

\vs 0.3cm \n We note that numerous properties of solutions to problem $(P)$, 
established in the two-dimensional case in \cite{[C]}, readily extend to higher dimension. In 
the ensuing discussion, we revisit some of these properties essential for our subsequent analysis.
The first property is the monotonicity of $\chi$ with respect to $x_n$.

\begin{proposition}\label{p1.2}
\begin{equation*}
\chi_{x_n}\leq 0 \quad\text{ in }
\quad{\cal D}^\prime(\Omega)
\end{equation*}
\end{proposition}

\vs 0,5cm \n Proposition 1.2, the continuity of $u$, and the maximum principle imply that if $u$ is positive at a 
given point of $\Omega$, then it remains positive throughout a cylinder below that point. 
Consequently, if $u$ vanishes at a specific point within $\Omega$, it also vanishes above that point.

\begin{proposition} \label{p1.3.} Let $x_0=(x_0',x_{0n})\in \Omega$.

\vs 0,2cm \n $i)\quad \text{If} \quad u(x_0) > 0,\quad \text{then there exists } \epsilon >0 \text{ such that}$
$$u(x',x_n) > 0 \qquad \forall
 (x',x_n)\in B_\epsilon(x_0)\cup \big\{ (x',x_n)\in \Omega\,/ \, \vert
 x'-x_0'\vert <\epsilon , \, x_n < x_{0n}  \big\},$$

\vs 0,2cm \n $ii)\quad \text{If} \quad u(x_0)=0,\quad  \text{then }
\quad u(x_0', x_n) =0 \qquad \forall x_n\geq x_{0n}.$
\end{proposition}

\vs 0,5cm \n Proposition 1.3 prompts the definition of the subsequent function for all $x'\in B_\rho'$:

\begin{equation*}
\phi(x') =\begin{cases}
& 0\qquad \text{ if }\quad \big\{x_n\in(0,l)\,/\, u(x',x_n)> 0~\big\}=\emptyset\\
& \sup\big\{x_n\in(0,l)\,/\, u(x',x_n)> 0~\big\}\qquad
\text{ otherwise}
\end{cases}
\end{equation*}

\vs 0.3cm\n $\phi$ is well defined and satisfies thanks to Proposition 1.3:

\vs 0.3cm
\begin{proposition}\label{p1.4.} $\phi$ is lower semi-continuous on $B_\rho'$ and
\[\big\{u>0\}=\{x_n<\phi(x')\big\}\]
\end{proposition}

\vs 0,5cm \n The last property states that if $u$ is zero in a particular open ball 
within $\Omega$, then $\chi$ is a.e. zero within that ball and above it. 

\begin{proposition}\label{p1.5}
Let $x_0=(x'_0,x_{0n})\in\Omega$ and $r>0$ such that
$B_r(x_0)\subset\subset \Omega $. If $u=0$ in $B_r(x_0)$, then
we have
\[\chi \, = 0 ~~\hbox{ a.e. in }
~~\Big\{\,(x',x_n)\in\Omega \,\,/\,\,\vert x' - x'_0
\vert < r\quad \hbox{and}\quad x_{0n} < x_n\,\Big\}\cup B_r(x_0)\]
\end{proposition}

\vs0.3cm \n The subsequent sections of this paper are structured as follows: 
In Section 2, we develop a super-solution for problem $(P)$, referred to as a "barrier function." 
This tool will play a pivotal role in proving Theorem 1.1.
Section 3 presents the detailed proofs of Theorem 1.1 and Corollary 1.1.

\section{A Barrier Function}\label{2}

\vs 0,5cm \n Let $x_0=(x'_0,x_{0n})\in \Omega$, $\overline{x}_0=(x_0',l)$, and $r_0>0$ such that
$\overline{B}_{2r_0}(x_0)\subset \Omega$, which implicitly implies that 
$\displaystyle{2r_0< \min(x_{0n},l-x_{0n})}$, 
and $K_0=\overline{B'}_{r_0}(x'_0)\times[x_{0n}-r_0,l]\subset\Omega\cup\Gamma$.

\vs 0,5cm \n From Proposition 1.1, there exists a constant $C_1>0$ depending on $n, p, \lambda, \Lambda, h, a, r_0, l$ such that:
\begin{equation}\label{e2.1}
|\nabla u(x)|\leq C_1\quad \text{ a.e. } x\in K_0
\end{equation}

\n We consider a function $\theta\in C^2(\mathbb{R})$ such that $-1\leq \theta\leq 1$ in $\mathbb{R}$ and
$$\theta(t)=\left\{
\begin{array}{ll}
-1 & \hbox{if } ~~t\leq \frac{1}{4}\\
1 & \hbox{if } ~~t\geq 1
    \end{array}
  \right.
$$
Setting $g(y')=\theta(|y'|^2)$ for $y'\in\mathbb{R}^{n-1}$, we have
$g\in C^2(\mathbb{R}^{n-1})$ and
$$g(y')=\left\{
    \begin{array}{ll}
      -1, & \hbox{if } ~~|y'|\leq \frac{1}{2}\\
      1, & \hbox{if } ~~|y'|\geq 1
    \end{array}
  \right.
$$

\n Let $r\in(0,r_0)$. Then the function $f_r(x')=l+rg\big({{x'-x'_0}\over r}\big)$ satisfies
$f_r\in C^2(\mathbb{R}^{n-1})$, $l-r\leq f_r(x')\leq l+r$ in $\mathbb{R}^{n-1}$, and
$$f_r(x')=\left\{
    \begin{array}{ll}
      l-r, & \hbox{if } ~~|x'-x'_0|\leq r/2\\
      l+r, & \hbox{if } ~~|x'-x'_0|\geq r.
    \end{array}
  \right.
$$

\n We define $\widetilde{\Omega}=B'_\rho\times(0,2l)$ and the extension 
of a function $f$ defined in $\Omega$ to $\widetilde{\Omega}$ by:
\[
\widetilde{f}(x',x_n) =\left\{
\begin{array}{ll}
f(x',x_n),  &\hbox{ if } (x',x_n)\in \Omega\\
f(x',2l-x_n),  &\hbox{ if } (x',x_n)\in \widetilde{\Omega}\setminus \Omega
       \end{array}
     \right.
\]

\n It is clear that $\widetilde{a}$ and $\widetilde{h}$ satisfy (1.1)-(1.5) and (1.6) in $\widetilde{\Omega}$
with the same constants.

\vs 0.1cm\n Finally, let $\epsilon\in(0,1)$ be a fixed real number. We consider the unique solution $v$ 
of the following elliptic Dirichlet boundary value problem
\begin{equation}\label{e2.2}
\left\{\begin{array}{ll}
div (\widetilde{a}(x) \nabla v +\widetilde{h} e_n)=0 & \hbox{in }~~ D_r \\
v={\epsilon\over 3}(f_r(x')-x_n) & \hbox{on }~~\partial D_r
    \end{array}
  \right.
\end{equation}

\n where $\displaystyle{ D_r=\{(x',x_n)\in\mathbb{R}^{n}~:~|x'-x'_0|<2r\text{ and }l-2r<x_n<f_r(x')~\}}$.

\n We observe that
\begin{equation}\label{e2.3}
0\leq v\leq{\epsilon\over 3}(l+r-(l-2r))=\epsilon r  ~~\hbox{ on }~\partial D_r.
\end{equation}

\vs0.3cm\n Let $G_r= \{x_n=f_r(x')\}\cap\{|x'-x'_0|\leq r\}$. Then we have:
\begin{proposition}\label{p2.1}
$v\in C^{1, \alpha}(D_r\cup G_r)$ and there exists a positive constant $C_2$  independent of $\epsilon$ and $r$ such that
for all $r\in(0,\min(r_0,\epsilon))$
$$ |\nabla v(x)|\,\leq\, C_2\epsilon^{1-{n\over p}}\qquad \forall x\in G_r$$

\end{proposition}

\vs0,3cm\n\emph{Proof.} Let $\displaystyle{\Delta=\Big\{(y',y_n)~:~|y'|<2\text{ and }-2<y_n<g(y')~\Big\}}$
and $r\in(0,\min(r_0,\epsilon))$. We introduce the
following functions defined on $\Delta$: 
\[w(y)={1\over r}v(\overline{x}_0+ry),\quad b(y)=\widetilde{a}(\overline{x}_0+ry),\quad  k(y)=\widetilde{h}(\overline{x}_0+ry),
\quad y={1\over r}(x-\overline{x}_0)
\]

\n Then it is straightforward to verify using equation (2.2) that $w$ satisfies
\begin{equation}\label{e2.4}
\left\{
\begin{array}{ll}
div(b(y) \nabla w)=- k_{y_n}=-r\widetilde{h}_{x_n}(\overline{x}_0+ry) & \hbox{in }~~ \Delta \\
w={\epsilon\over 3}(g(y')-y_n) & \hbox{on }~~\partial\Delta.
    \end{array}
  \right.
\end{equation}
Clearly $b$ and $k$ satisfy (1.1)-(1.6) with the same constants.
Moreover, we have from (1.5) and (2.3)-(2.4)
\begin{equation*}
\left\{
\begin{array}{ll}
div (b(y) \nabla w)\leq 0 & \hbox{in }~~ \Delta \\
0\leq w\leq \epsilon & \hbox{on }~~ \partial \Delta
    \end{array}
  \right.
\end{equation*}
which in turn leads, using the weak maximum principle, to
\begin{equation}\label{e2.5}
0\leq w\leq \epsilon~\hbox{ in }~ \Delta
\end{equation}

\vs 0.2cm\n Let $T=\big\{y_n=g(y')\big\}\cap\big\{|y'|\leq 1\big\}$ and
$\Delta'=\Delta\cap\big\{|y'|\leq 3/2\big\}$. Since $T$ is a $C^{1, \alpha}$
boundary portion of $\partial \Delta$, $w=0$ on $T$, we deduce from (2.4), and
the subsequent Remark to Corollary 8.36 p. 212, reference \cite{[GT]} that
$w\in C^{1, \alpha}(\Delta\cup T)$ and satisfies the following estimate

\begin{equation}\label{e2.6}
|w|_{1,\alpha, \Delta'}\leq C \Big(|w|_{0,\Delta}
  \,+\, |k_{y_n}|_{p,\Delta}\Big)
\end{equation}

\n Moreover, we know from the above reference that $C=C(n,\lambda,M,d',T)$, $~d'=dist(\Delta',\partial\Delta\setminus T)$,
$M$ is an upper bound of
\begin{eqnarray*}
&&\displaystyle{\max_{i,j}(|b_{ij}|_{0,\alpha, \Delta'})=
r^\alpha\max_{i,j}(|a_{ij}|_{0,\alpha, D_r'})\leq \max_{i,j}(|\widetilde{a}_{ij}|_{0,\alpha, \widetilde{K}_0})}, \quad\text{since } r<\epsilon<1\\
&&D_r'=D_r\cap\{l-3r/2\leq x_n\leq l+r/2\}\subset\subset \widetilde{K}_0=\overline{B}_{r_0}'(x_0')\times[x_{0n}-r_0,2l]
\end{eqnarray*} 

\n Since $\displaystyle{\max_{i,j}(|\widetilde{a}_{ij}|_{0,\alpha, \widetilde{K}_0})=\max_{i,j}(|a_{ij}|_{0,\alpha, K_0})}$, 
it follows that $C$ can be chosen to depend solely on $n$, $\lambda$, and $\displaystyle{\max_{i,j}(|a_{ij}|_{0,\alpha, K_0})}$.

\vs0.1cm\n Now we have
\begin{eqnarray}\label{e2.7}
|k_{y_n}|_{p,\Delta}&=&|r\widetilde{h}_{x_n}(\overline{x}_0+ry)|_{p,\Delta}=\Bigl(\int_\Delta
r^p \widetilde{h}_{x_n}^p(\overline{x}_0+ry)dy\Bigl)^{1/p}\nonumber\\
&=&\Bigl(\int_{D_r}{{r^p}\over {r^n}}\widetilde{h}_{x_n}^p(x)dx\Bigl)^{1/p}=r^{1-{n\over p}}|\widetilde{h}_{x_n}|_{p,D_r}
\end{eqnarray}

\n Hence, we infer from (2.5)-(2.7) and the fact that $\epsilon<1$, that we have for all $0<r<\min(\epsilon,r_0)$
\begin{eqnarray*}
|w|_{1,\alpha, \Delta'}\leq C \Big(\epsilon+\epsilon^{1-{n\over p}}|\widetilde{h}_{x_n}|_{p,\widetilde{\Omega}}\Big)
\leq C \Big(\epsilon^{{n\over p}}+|\widetilde{h}_{x_n}|_{p,\widetilde{\Omega}}\Big)\epsilon^{1-{n\over p}}
\leq C \Big(1+|\widetilde{h}_{x_n}|_{p,\widetilde{\Omega}}\Big)\epsilon^{1-{n\over p}} 
\end{eqnarray*}
\n We obtain for another constant $\displaystyle{C=C(n, \lambda, |h_{x_n}|_{p,\Omega}, \max_{i,j}(|a_{ij}|_{0,\alpha, K_0}))}$ 
denoted by $C_2$, $~|\nabla w|_{0,\Delta'}\leq C_2\epsilon^{1-{n\over p}}$. This leads, in particular, to
$$|\nabla w(y',g(y'))|\leq C_2\epsilon^{1-{n\over p}}\qquad \forall y'\in \overline{B}'_1$$
Therefore, we obtain 
\[|\nabla v(x',f_r(x'))|=\big|\nabla
w\big({{x'-x'_0}\over r},g\big({{x'-x'_0}\over r}\big)\big |
\,\leq\, C_2\, \epsilon^{1-{n\over p}}\qquad \forall x'\in G_r\] \qed

\vs 0,2cm\n Let now $\Omega_r=\{(x',x_n)\in\Omega~:~|x'-x'_0|<r \text{ and } l-2r<x_n<l~\}$,
$\displaystyle{\theta_0=2\max_{0\leq t\leq 1}|\theta'(t)|}$ and $\epsilon_0=\Big({{\underline{h}}\over{C \Lambda\sqrt{1+\theta_0^2}}}\Big)^{p\over{p-n}}$.  

\vs 0,2cm\n From now on, we will denote by $v$ the function defined in $(2.1)$ for $\epsilon=\epsilon_0$
and extended by $0$ to $\big\{x_n>f_r(x')\big\}$. The next lemma asserts that $v$ is a super-solution to problem $(P)$ in $\Omega_r$.

\begin{lemma}\label{l2.1} We have for all $0<r<\min(r_0,\epsilon_0)=r_1$:
\begin{eqnarray}\label{e2.8}
&&\int_{\Omega_r}\Big(a(x)  \nabla v + \chi_{\{v>0\}}h(x)
e_n\Big).\nabla\zeta dx\geq 0 \nonumber\\
&&\hskip 1.5cm\forall \zeta\in H^1(\Omega_r),\quad \zeta\geq 0 \mbox{ in }
\Omega_r,\quad\zeta=0 \mbox{ on }
\partial \Omega_r\setminus\Gamma
\end{eqnarray}
\end{lemma}

\vs 0,3cm \n\emph{Proof.} Set $C_r=\Omega_r\cap\{x_n<f_r(x')\}=\Omega_r\cap\{v>0\}$,
and let $\nu=(\nu_1,...,\nu_n)$ be the outward unit normal vector to $\partial C_r\cap\left(\{x_n=l\}\cup\Omega_r\right)$.
First, we have for $|x'-x'_0|\leq r$
\begin{equation*}
|\nabla_{x'}f_r(x')|=\Big|\theta'\Big({{|x'-x'_0|^2}\over r^2}\Big)\Big|.{{2|x'-x'_0|}\over r}\leq \theta_0,
\end{equation*}
which leads to
\begin{equation}\label{e2.9}
\nu_{n}={1\over{\sqrt{1+|\nabla_{x'}f_r(x')|^2}}}\geq {1\over{\sqrt{1+\theta_0^2}}},~\text{ on }~\partial C_r\cap\Omega_r=\{x_n=f_r(x')\}\cap\Omega_r
\end{equation}

\n Next, since $\displaystyle{\nu_n=1\geq {1\over{\sqrt{1+\theta_0^2}}}}$ on $\{x_n=l\}\cap\partial C_r$, we deduce from (1.1), (1.3), Proposition 2.1, and (2.9), that for $0<r<r_1$
\begin{eqnarray}\label{e2.10}
 a(x) \nabla v. \nu +h(x)\nu_{n}\geq-C_1 \Lambda.\epsilon_0^{1-{n\over p}} +{{\underline{h}}\over{\sqrt{1+\theta_0^2}}}=0,
~\text{ on }~\partial C_r\cap\left(\{x_n=l\}\cup\Omega_r\right)
\end{eqnarray}

\n We conclude by taking into account (2.1) and (2.10), for $\zeta\in H^1(\Omega_r)$ such that $\zeta\geq 0$ in 
$\Omega_r$, $\zeta=0$ on $\partial \Omega_r\setminus\Gamma$

\begin{eqnarray*} && \int_{\Omega_r}\Big(a(x)  \nabla v +\chi_{\{v>0\}}h(x)
 e_n\Big).\nabla\zeta dx \\
&&=\int_{C_r}\Big(a(x) \nabla v +h(x)e_n\Big).\nabla\zeta dx\\
&&=\int_{\partial C_r\cap\left(\{x_n=l\}\cup\Omega_r\right)}\Big(a(x) \nabla v.\nu + h(x)
\nu_n\Big)\zeta d\sigma \geq 0
\end{eqnarray*}\qed

\vs 0,3cm
\begin{remark}\label{r2.1}

\n $(i)$ Functions like the function $v$ constructed earlier are commonly referred to as a barrier functions.

\vs 0.2cm \n $(ii)$ The proof of Theorem 1.1 is based on comparing the function
$u$ to the barrier function $v$. We initiate this process by establishing two lemmas. 
\end{remark}

\vs 0,3cm
\begin{lemma}\label{l2.2}
Assume that $ u\leq v$ on $\Omega\cap\partial \Omega_r$, and let $\Omega_r^\delta=\Omega_r \cap \{0\leq u-v< \delta\}\cap \{v>0\}$.
Then we have
\begin{equation}\label{e2.11}
\int_{\Omega_r^\delta } a(x)  \nabla (u-v)^+. \nabla (u- v)^+ dx \,=\, o(\delta)~~{as }~~\delta\rightarrow 0,
\quad \text{uniformly in } r\in(0,r_1]
\end{equation}

\end{lemma}

\vs 0,3cm\n\emph{Proof.} For $\delta,\eta> 0$, we consider the functions
$$H_\delta (s)=\displaystyle{\min\Big({s^+\over \delta},1\Big)},\qquad d_\eta (x)=H_\eta ((x_n-f_r(x'))^+).$$

\n Given that  $u\leq v$ on $\partial \Omega_r\cap\Omega$, and since $u=0$ on $\partial \Omega_r\cap\Gamma$
and $H_\delta(s)=0$ for $s\leq 0$, we have $H_\delta(u-v)=0$ on
$\partial \Omega_r$. Moreover, since $d_\eta=0$ below $\{x_n=f_r(x')\}$, we
have $d_\eta=0$ on $\Omega\cap\partial \Omega_r$. Therefore
$\zeta= H_\delta (u-v) + d_\eta  (1- H_\delta (u))\in H^1(\Omega_r)$ is a
nonnegative function that vanishes on $\Omega\cap\partial \Omega_r$, hence we have
\begin{eqnarray*}
\int_{\Omega_r}\Big(a(x)  \nabla u + \chi h(x)e_n\Big).\nabla \zeta dx \leq 0,
\end{eqnarray*}
which leads to
\begin{eqnarray}\label{e2.12}
&&\int_{\Omega_r}\Big(a(x)  \nabla u + \chi h(x)
 e_n\Big).\nabla (H_\delta (u-v) ) dx \nonumber\\
&&\qquad \leq - \int_{\Omega_r}\Big(a(x)  \nabla u +\chi h(x)
 e_n\Big).\nabla (d_\eta  (1- H_\delta (u)) ) dx
\end{eqnarray}

\n Since $H_\delta(u-v)=0$ on $\partial \Omega_r$, we get by (2.8)
\begin{equation}\label{e2.13}
-\int_{\Omega_r}\Big(a(x)  \nabla v +\chi_{\{v>0\}}h(x)
e_n\Big).\nabla (H_\delta (u-v) ) dx\leq 0
\end{equation}
Adding (2.12) and (2.13), we get
\begin{eqnarray*}
&&\int_{\Omega_r}\Big(a(x)\nabla (u-v) + (\chi-\chi_{\{v>0\}}) h(x)
 e_n\Big).\nabla (H_\delta (u-v) ) dx \nonumber\\
&&\qquad \leq - \int_{\Omega_r}\Big(a(x)  \nabla u +\chi h(x)
 e_n\Big).\nabla (d_\eta  (1- H_\delta (u)) ) dx
\end{eqnarray*}
which can be written as
\begin{eqnarray*}
 &&\int_{\Omega_r}  a(x)  \nabla (u-v). \nabla (H_\delta (u-v) )
dx \leq \int_{\Omega_r} h(x) (\chi_{\{v>0\}}-\chi) e_n. \nabla (H_\delta(u-v)) dx\nonumber\\
&&\qquad\qquad -\int_{\Omega_r}\Big(a(x)  \nabla u + \chi h(x)
 e_n\Big).\nabla (d_\eta  (1- H_\delta (u)) ) dx
\end{eqnarray*}
or
\begin{eqnarray}\label{e2.14}
 &&\int_{\Omega_r\cap \{v>0\}}  a(x)  \nabla (u-v). \nabla (H_\delta (u-v) )
dx \leq \int_{\Omega_r} h(x) (\chi_{\{v>0\}}-\chi) e_n. \nabla (H_\delta(u-v)) dx\nonumber\\
&&\quad-\int_{\Omega_r\cap \{v=0\}}  a(x)  \nabla u. \nabla (H_\delta (u) )
dx -\int_{\Omega_r}\Big(a(x)  \nabla u + \chi h(x)e_n\Big).\nabla (d_\eta  (1- H_\delta (u)) ) dx\nonumber\\
&&
\end{eqnarray}
Since $H_\delta (u-v)=0$ in $\{u\leq v\}$, and $\chi=1$ a.e. in $\{u>0\}\supset\{u>v\}$,
we have
\begin{eqnarray}\label{e2.15}
&&\int_{\Omega_r}(\chi_{\{v>0\}}-\chi) h(x)e_n.\nabla (H_\delta (u-v) ) dx=-\int_{\Omega_r\cap \{v=0\}} \chi h(x)e_n.\nabla (H_\delta (u) ) dx\nonumber\\
&&\hskip 2cm=\int_{\Omega_r\cap \{v=0\}} \chi h(x)e_n.\nabla (1- H_\delta (u))dx
\end{eqnarray}

\n Given that $d_\eta=0$ in $\{v>0\}$, we have
\begin{eqnarray}\label{e2.16}
&&\int_{\Omega_r}\Big(a(x)  \nabla u + \chi h(x)e_n\Big).\nabla (d_\eta  (1- H_\delta (u)) ) dx\nonumber\\
&&\hskip 2cm=\int_{\Omega_r\cap \{v=0\}}\Big(a(x)  \nabla u + \chi h(x)e_n\Big).\nabla (d_\eta  (1- H_\delta (u)) ) dx
\end{eqnarray}

\n Taking into account (2.14)-(2.16) and the fact that $\nabla (1-H_\delta (u))=-\nabla (H_\delta(u))$, we obtain
\begin{eqnarray} \label{e2.17}
&&\int_{\Omega_r\cap \{v>0\}} H_\delta'(u-v)  a(x)  \nabla (u-v). \nabla (u-v) dx
\leq \int_{\Omega_r\cap \{v=0\}}  a(x)  \nabla u. \nabla (1- H_\delta (u))dx\nonumber\\
&&\quad+\int_{\Omega_r\cap \{v=0\}} \chi h(x)e_n.\nabla (1- H_\delta (u))dx\nonumber\\
&&-\int_{\Omega_r\cap \{v=0\}}\Big(a(x)  \nabla u + \chi h(x)e_n\Big).\nabla (d_\eta  (1- H_\delta (u)) ) dx\nonumber\\
&&=\int_{\Omega_r\cap \{v=0\}}\Big(a(x)  \nabla u + \chi h(x)
 e_n\Big).\nabla ((1-d_\eta)  (1- H_\delta (u)) ) dx\nonumber\\
&&\quad= -\int_{\Omega_r\cap \{v=0\}} (1- d_\eta)\Big(a(x)  \nabla u +\chi h(x)
 e_n\Big).\nabla ( H_\delta (u) ) dx\nonumber\\
&&\quad - \int_{\Omega_r\cap \{v=0\}} (1- H_\delta (u))\Big(a(x)  \nabla u +\chi h(x)
 e_n\Big).\nabla d_\eta dx  = I_1^{\delta,\eta} + I_2^{\delta,\eta}
\end{eqnarray}
Since $\displaystyle{\lim_{\eta\rightarrow 0}d_\eta(x)=1}$ for a.e. $x\in \Omega_r\cap \{v=0\}$,
we have by the Lebesgue theorem
\begin{equation}\label{e2.18}
\lim_{\eta\rightarrow 0} I_1^{\delta,\eta}=0
\end{equation}
\n Now, observe that we have $\nabla d_\eta=H'_\eta((x_n-f_r(x'))^+)(-\nabla_{x'}f_r(x'),1)$
in $\Omega_r\cap\{v=0\}$. In particular, we see that $d_\eta$ is an increasing function in $x_n$, leading to
\begin{eqnarray}\label{e2.19}
I_2^{\delta,\eta}&=& -\int_{\Omega_r\cap \{u=v=0\}} \chi h(x)
 e_n.\nabla d_\eta dx\nonumber\\
&&- \int_{\Omega_r\cap \{u>v=0\}} (1- H_\delta (u))\big(a(x)
\nabla u + h(x)e_n\big).\nabla d_\eta dx\nonumber\\
&=&-\int_{\Omega_r\cap \{u=v=0\}} H'_\eta((x_n-f_r(x'))^+)\chi h(x) dx\nonumber\\
&&- \int_{\Omega_r\cap \{u>v=0\}} (1- H_\delta (u))\big(a(x)
\nabla u + h(x)e_n\big).\nabla d_\eta dx\nonumber\\
&\leq& - \int_{\Omega_r\cap \{u>v=0\}} (1- H_\delta (u))\big(a(x)
\nabla u + h(x)e_n\big).\nabla d_\eta dx
\end{eqnarray}
\n Using (1.3) and (2.1), we get from (2.19), for $C_3=\Lambda C_1+\overline{h}$
\begin{eqnarray}\label{e2.20}
| I_2^{\delta,\eta}|&\leq& \displaystyle{{{C_3}\over \eta}
\int_{\{u>0\}\cap \{f_r(x')< x_n< f_r(x')+\eta\}}
 (1- H_\delta (u)) }dx\nonumber\\
 &\leq& \displaystyle{{C_3\over \eta} \int_{\{f_r(x')< x_n< f_r(x')+\eta\}\cap\{f_r(x')<\phi(x')\}}
(1- H_\delta (u))}dx\nonumber\\
 &=& C_3 \int_{\{f_r<\phi\}}\Big({1\over \eta}\int_{f_r(x')}^{f_r(x')+\eta}(1- H_\delta (u))dx_n\Big)dx'
\end{eqnarray}
\n For each $x'\in \{f_r<\phi\}$, we have by the continuity of $1- H_\delta (u)$
\begin{equation}\label{e2.21}
\lim_{\eta\rightarrow 0}{1\over \eta}\int_{f_r(x')}^{f_r(x')+\eta}(1- H_\delta (u))dx_n=1- H_\delta (u(x',f_r(x')))
\end{equation}
Using (2.20)-(2.21) and the Lebesgue theorem, we obtain by letting $\eta\rightarrow0$
\begin{equation}\label{e2.22}
\limsup_{\eta\rightarrow 0}| I_2^{\delta,\eta} | \leq
C_3\int_{\{f_r<\phi\}}(1- H_\delta (u(x',f_r(x'))))dx'
\end{equation}
Hence, combining (2.17)-(2.18) and (2.22), we arrive at
\begin{eqnarray*}
&&{1\over\delta }\int_{\Omega_r\cap\{0<u-v< \delta\}\cap \{v>0\}}
 a(x)  \nabla (u-v)^+. \nabla (u-v )^+dx\\
&&\quad\leq  C_3\int_{\{f_r<\phi\}}(1- H_\delta (u(x',f_r(x'))))dx'
\end{eqnarray*}

\n Since $\{f_r<\phi\}\subset \{u>0\}$, we have
$\displaystyle{\lim_{\delta\rightarrow 0}}
(1- H_\delta (u(x',f_r(x'))))=0$. The lemma follows by using the Lebesgue 
theorem in the above inequality.\qed

\vs 0,5cm
\begin{lemma}\label{l.2.3} If $u\leq v$ on $\partial \Omega_r\cap\Omega$, then
$u\equiv 0$ in $C_{r/2}(x_0)\cap \{x_n\geq l-r\}$.
\end{lemma}

\vs 0,3cm \n\emph{Proof.} Let $\displaystyle{\zeta \in \mathcal{D}(\Omega\cap\{x_n<f^*_r(x')\})}$, where $f^*_r(x')=\min(l,f_r(x')$.
Then we have by Lebesgue's theorem
\begin{eqnarray}\label{e2.23}
\int_{\Omega_r} a(x) \nabla (u-v)^+\cdot\nabla\zeta dx
&=& \,\lim_{\delta\rightarrow 0} \int_{\Omega_r}H_\delta (u-v) a(x)
\nabla (u-v)^+\cdot\nabla\zeta dx\nonumber\\
&=& \,\lim_{\delta\rightarrow 0}I_\delta
\end{eqnarray}
\n Since $H_\delta (u-v)=0$ whenever $u\leq v$, we have
$\nabla(u-v)^+\cdot\nabla (H_\delta (u-v)\zeta)=\nabla (u-v)\cdot\nabla
(H_\delta (u-v)\zeta)$ a.e. in $\Omega_r$. Therefore, we can write
\begin{eqnarray*}
I_\delta &=&  \int_{\Omega_r} a(x) \nabla (u-v)\cdot\nabla (H_\delta (u-v)\zeta)dx\\
&&-\int_{\Omega_r} \zeta a(x) \nabla (u-v)\cdot\nabla(H_\delta (u-v))dx= I^1_\delta - I^2_\delta
\end{eqnarray*}

\n or
\begin{eqnarray}\label{e2.24}
I_\delta &=&  \int_{\Omega_r} a(x) \nabla(u-v)\cdot\nabla (H_\delta (u-v)\zeta)dx\nonumber\\
&&\quad-{1\over \delta}
\int_{\Omega_r \cap \{0<u-v<\delta\}\cap\{v>0\}} \zeta a(x) \nabla (u-v)\cdot\nabla  (u-v)dx\nonumber\\
&=& I^1_\delta - I^2_\delta
\end{eqnarray}

\n First, we observe by using (1.2), that
$$|I^2_\delta | \leq \max_{\overline{\Omega_r}} |\zeta| \, . \,{1\over \delta}
\int_{\Omega_r \cap \{0<u-v<\delta\}\cap\{v>0\}} a(x) \nabla (u-v)\cdot\nabla(u-v) dx$$

\n which leads by Lemma 2.2 to

\begin{equation}\label{e2.25}
\lim_{\delta\rightarrow 0} I_\delta^2=0
\end{equation}

\n Next, we claim that $I^1_\delta=0$. Indeed,

\begin{equation}\label{e2.26}
I^1_\delta= \int_{\Omega_r} a(x) \nabla u\cdot\nabla (H_\delta
(u-v)\zeta)dx - \int_{\Omega_r} a(x) \nabla v\cdot\nabla (H_\delta
(u-v)\zeta)dx.
\end{equation}

\n Since $u\leq v$ on $\Omega\cap\partial \Omega_r$, we have $H_\delta(u-v)=0$ on
$\Omega\cap \partial \Omega_r$. Moreover, $\displaystyle{\zeta \in \mathcal{D}(\Omega\cap\{x_n<f^*_r(x')\})}$
and $\displaystyle{\Omega_r\cap\{v>0\}=\Omega_r\cap\{x_n<f^*_r(x')\} }$ imply that
$H_\delta(u-v)\zeta\in H^1_0(\Omega_r\cap\{x_n<f^*_r(x')\})$. Therefore, we obtain from (2.2)
\begin{eqnarray}\label{e2.27}
\int_{\Omega_r} a(x) \nabla v\cdot\nabla (H_\delta (u-v)\zeta)dx
&=&\int_{\Omega_r\cap\{x_n<f^*_r(x')\}} a(x) \nabla v\cdot\nabla (H_\delta (u-v)\zeta)dx\nonumber\\
&=&-\int_{\Omega_r\cap\{x_n<f^*_r(x')\}} h(x)(H_\delta (u-v)\zeta)_{x_n} dx\nonumber\\
&=&-\int_{\Omega_r} h(x)(H_\delta (u-v)\zeta)_{x_n} dx
\end{eqnarray}
\n Now, since $\pm \chi(\Omega_r) H_\delta(u-v)\zeta$ are test functions for $(P)$, $\chi=1$ a.e. in 
$\Omega_r\cap\{u>0\}$ and $H_\delta (u-v)=0$ whenever $u=0$, we have

\begin{eqnarray}\label{e3.28}
&&\int_{\Omega_r} a(x) \nabla u. \nabla (H_\delta(u-v)\zeta)dx =-\int_{\Omega_r} \chi h(x)(H_\delta
(u-v)\zeta)_{x_n}dx\nonumber\\
&&\quad=-\int_{\Omega_r}h(x)(H_\delta(u-v)\zeta)_{x_n} dx
\end{eqnarray}

\n It follows from (2.26)-(2.28) that $I^1_\delta=0$. Hence, we get from
(2.23)-(2.25) that

$$\int_{\Omega_r} a(x) \nabla (u-v)^+\cdot\nabla\zeta dx\, = \,0 \qquad
\forall\zeta \in  \mathcal{D}(\Omega\cap\{x_n<f^*_r(x')\})$$

\n Given that $(u-v)^+=0$ on $\partial \Omega_r$, the extension $w=\chi(\Omega_r)(u-v)^+$ by 0
of $(u-v)^+$ to $\Omega\setminus \Omega_r$ belongs to $H^1(\Omega)$ and we obtain
\[\int_{\Omega\cap\{x_n<f^*_r(x')\}} a(x) \nabla w. \nabla\zeta dx\, = \,0 \qquad
\forall\zeta \in  \mathcal{D}(\Omega\cap\{x_n<f^*_r(x')\})\]
Taking into account (1.1)-(1.2), and the fact that $w\geq 0$ in $\Omega$ and
$w=0$ in $\Omega\setminus \Omega_r$, we get by the strong maximum principle that
$w=0$ in $\Omega\cap\{x_n<f^*_r(x')\}$,
which leads to $(u-v)^+=0$ in $\displaystyle{\Omega_r\cap\{x_n<f^*_r(x')\}}$,
or equivalently $u\leq v$ in $\displaystyle{\Omega_r\cap\{x_n<f^*_r(x')\}}$.
We recall that $v=0$ on $\Omega_r\cap\{x_n=f^*_r(x')\}$. Therefore, we get
$u=0$ on $\Omega_r\cap\{x_n=f^*_r(x')\}$, and consequently, we obtain by using
Proposition 1.3 $(ii)$ that $u=0$ in $\displaystyle{\Omega_r\cap\{x_n\geq f^*_r(x')\}}$.
Since $f^*_r(x')=l-r$ when $|x'-x'_0|<r/2$, we conclude that 
$u=0$ in the cylinder $C_{r/2}\cap \{x_n\geq l-r\}$.
\qed

\vs 0,5cm

\section{Proof of the Main Results}\label{3}

\vs 0,3cm \n\emph{Proof of Theorem 1.1.} Let $x_0=(x'_0,x_{0n})\in FB$ and let $r\in(0,r_2)$. 
We continue to use the same notation as in the previous section. We know that
\begin{equation}\label{3.1}
v_{/\partial \Omega_r\cap \Omega}={\epsilon_0\over 3}(f_r(x')-x_n)\geq{\epsilon_0\over 3}(l+r-l)=\frac{\epsilon_0 r}{3}=C_0r
\end{equation}

\n According to Proposition 1.3 ($i$), we have $u(x'_0,x_n)>0$ for all $x_n \in (0,x_{0n})$. 
Specifically, $u$ is not identically zero within $C_{r/2}\cap \{x_n\geq l-r\}$.
Consequently, from Lemma 2.3 and (3.1), we infer that $\displaystyle{\sup_{\overline{\Omega}_r}u > C_0r}$. 
Thus, there exists $x_1=(x'_1,x_{1n}) \in \overline{\Omega}_r$ such that $u(x_1) > C_0r$.
Using (2.1), we further conclude that for each $x$ in $B_{\delta r}(x_1)$, where $\displaystyle{\delta=\min\left(1,\frac{C_0}{2C_1}\right)}$
\begin{eqnarray*}
u(x)\geq u(x_1)-C_1|x-x_1|\geq C_0r-C_1\delta r\geq\Big(C_0-C_1{C_0\over{2C_1}}\Big) r={C_0\over2}r
\end{eqnarray*}

\n We deduce that $B_{\delta r}(x_1)$ is contained within the set $\{u > 0\}$, which implies by Proposition 1.3 ($i$)
that $u>0$ below $B_{\delta r}(x_1)$. Picking the point $x_2=(x'_1,x_{0n})$, it follows that $B_{\delta r}(x_2)$ also 
lies within $\{u > 0\}$. Furthermore, since $\delta\leq 1$, $|x'_2-x_0'|=|x'_1-x_0'|\leq r$, and $x_{2n}=x_{0n}$, it is
easy to verify that $B_{\delta r}(x_2) \subset B_{2r}(x_0)$. Thus, we have established that 

\begin{equation}\label{3.2}
B_{\delta r}(x_2) \subset B_{2r}(x_0) \setminus FB\quad \forall r\in(0,r_2)
\end{equation}

\n To conclude, let $\delta_0=\delta/2$, and use (3.2) for $r/2$ with $r\in(0,r_2)$, we get since $r/2\in(0,r_2)$
\[
B_{\delta_0 r}(x_2)= B_{\delta r/2}(x_2) \subset B_{2(r/2)}(x_0) \setminus FB=B_{r}(x_0) \setminus FB\quad \forall r\in(0,r_2)\]

\n This completes the proof, demonstrating that the free boundary is a porous set with porosity $\delta_0$.
\qed

\vs 0,5cm \n\emph{Proof of Corollary 1.1.} First we have by $(P)(ii)$
\begin{equation}\label{e3.3}
\chi=1~~~ \text{a.e. in }~~\{u>0\}.
\end{equation}
\n Next, let $x$ be a point in the interior of the set $\{u=0\}$ denoted by $Int(\{u=0\})$. 
There exists a ball $B_r(x)$ of center $x$ and radius $r$ such that $u=0$ in $B_r(x)$.
By Proposition 1.5, we have $\chi=0$ a.e. in $B_r(x)$. Therefore, we obtain
\begin{equation}\label{e3.4}
\chi=0~~~ \text{a.e. in }~~Int(\{u=0\})
\end{equation}

\n Since $FB=\partial\{u>0\}\cap\Omega$ is of Lebesgue's measure zero due to its porosity, we conclude from (3.3)-(3.4)
that $\chi=\chi_{\{u>0\}}$ a.e. in   $\Omega$.
\qed

\vs 0.5cm\n 
\begin{remark}\label{r3.1} 

As stated in the introduction, when $h(x)=a_{nn}(x)$, the problem $(P)$ transforms 
into the dam problem. In \cite{[L1]}, the author established that the so-called Reservoirs-Connected solution 
of the dam problem is unique provided that the set $FB$ has Lebesgue measure zero for $n\geq 2$. 
Therefore, this uniqueness directly stems from Theorem 1.1.

\end{remark}

\end{document}